\input amstex
\documentstyle{amsppt}
\magnification=\magstep1
\hoffset=1.25truein
\voffset=2truein
\vsize=7.25truein
\hsize=5.377593truein


\overfullrule=0pt
\nologo

\def\dotS{{\dot S}}
\def\cL{{\Cal L}}
\def\cP{{\Cal P}}
\def\cV{{\Cal V}}
\def\oV{{\overline{V}}}
\def\oZ{{\overline{Z}}}

\topmatter
\title{On the hereditary paracompactness of locally compact,
hereditarily normal spaces}\endtitle
\author{Paul Larson$^1$ and
    Franklin D. Tall$^2$}
\endauthor

\date{November 29, 2010}\enddate

\rightheadtext{Paul Larson$^1$ and Franklin D. Tall$^2$}

\address{Paul Larson, Department of Mathematics,
University of Toronto, Toronto, Ontario  M5S 3G3 CANADA}
\endaddress
\email{larson\@math.utoronto.ca}\endemail
\address{Franklin D. Tall, Department of Mathematics,
University of Toronto, Toronto, Ontario  M5S 3G3 CANADA}
\endaddress
\email{f.tall\@utoronto.ca}\endemail

{\footnotetext""{AMS Subj. Class. (2010):
Primary 54D35, 54D15, 54D20, 54D45, 03E65;
Secondary 03E35.}}
{\footnotetext""{{\it Key words and phrases.} locally compact,
hereditarily normal, paracompact,
Axiom R, PFA$^{++}$.}}
{\footnotetext"$^1$"{\ The first author acknowledges support from
Centre de Recerca Mathem\`atica and from NSF-DMS-0801009.}}
{\footnotetext"$^2$"{\ The second author acknowledges support from
NSERC grant A-7354.}}

\abstract
We establish that if it is consistent that there is a
supercompact cardinal, then it is consistent that every locally
compact, hereditarily normal space which does not include a perfect
pre-image of $\omega_1$ is hereditarily paracompact.
\endabstract
\endtopmatter

\document

This is the fifth in a series of papers
(\cite{LTo}, \cite{L$_2$}, \cite{FTT}, \cite{LT}, \cite{T$_3$} being the logically previous
ones) that establish powerful topological consequences in models
of set theory obtained by starting with a particular kind of Souslin
tree $S$, iterating partial orders that don't destroy $S$, and then
forcing with $S$.  The particular case of the theorem stated in
the abstract when $X$ is perfectly normal (and hence has no perfect
pre-image of $\omega_1$) was proved in \cite{LT},
using essentially that locally compact perfectly normal spaces are
locally hereditarily Lindel\"of and first countable.
Here we avoid these two last properties by combining the methods of
\cite{B$_2$} and \cite{T$_3$}.
To apply \cite{B$_2$}, we establish the new set-theoretic result that
PFA$^{++}(S)[S]$ implies Fleissner's ``Axiom R''.
This notation is explained below; the model is a strengthening of those
used in the previous four papers.
\medskip

The results established here were actually proved around 2004, modulo
results of Todorcevic announced in 2002 (which now appear in \cite{FTT} and \cite{L$_2$}) and of the second author \cite{T$_3$}. We
have delayed submission until a correct version of \cite{T$_3$}
existed in preprint form.
\medskip

\proclaim{Definition}
A continuous map is {\bf perfect} if images of closed sets are closed, and pre-images of points are compact.
\endproclaim

It is easy to find locally compact, hereditarily normal spaces
which are not paracompact -- $\omega_1$ is one such.  Non-trivial
perfect pre-images of $\omega_1$ may also be hereditarily normal,
but are not paracompact.  Our result says that consistently, any
example must in fact include such a canonical example.
\medskip

\proclaim{Theorem 1}
If it is consistent that there is a supercompact cardinal,
it's consistent that every locally compact, hereditarily
normal space that does not include a perfect
pre-image of $\omega_1$ is (hereditarily) paracompact.
\endproclaim
\medskip

This is not a ZFC result, since there are many consistent
examples of locally compact, perfectly normal spaces which are
not paracompact.  For example, the Cantor tree over a $Q$-set, which is the standard example of a locally compact, normal, non-metrizable Moore space -- see e.g. \cite{T$_1$}, which has essentially the same example.  Other examples include the Ostaszewski and Kunen lines, as in \cite{FH}.
\medskip

Let us state some axioms we will be using.
\medskip

\noindent {\bf PFA$^{++}$:} {\it
Suppose $P$ is a proper partial order,
$\{ D_\alpha \}_{\alpha < \omega_1} $ is a collection
of dense subsets of $P$, and
$\{  \dotS_\alpha : \alpha <\omega_1 \}$ is a sequence of
terms such that $(\forall \alpha < \omega_1 ) \Vdash_P \dotS_\alpha$
is stationary in $\omega_1$.  Then there is a filter
$G\subseteq P$ such that
\roster
\item "{\it (i)}"
$(\forall \alpha < \omega_1 )\  G\cap D_\alpha \ne 0 $,
\item "{\it and (ii)}"
$(\forall \alpha < \omega_1 )\ S_\alpha (G) = \{ \xi < \omega_1 :
(\exists p \in G ) p \Vdash \xi \in \dotS_\alpha \} $
is stationary in $\omega_1$.
\endroster
}
\medskip

Baumgartner \cite{Ba} introduced this axiom and called it
``PFA$^{+}$".  Since then, others have called this
``PFA$^{++}$", using ``PFA$^{+}$" for the weaker one-term version.
As Baumgartner observed, the usual consistency proof for PFA,
which uses a supercompact cardinal, yields
a model for what we are calling PFA$^{++}$.
\medskip

\proclaim{Definition}
$\Gamma \subseteq [X]^{<\kappa}$ is {\bf tight} if whenever
$\{ C_\alpha : \alpha < \delta \}$ is an increasing
sequence from $\Gamma$,  and $\omega < cf \delta < \kappa $,
then $\bigcup \{ C_\alpha : \alpha < \delta \} \in \Gamma $.
{\bf Axiom R:} if $\Sigma \subseteq [X]^{<\omega_1} $ is stationary
and $\Gamma \subseteq [X]^{<\omega_2} $ is tight and cofinal,
then there is a $Y \in \Gamma $ such that $\cP(Y) \cap \Sigma $
is stationary in $[Y]^{<\omega_1} $.
{\bf Axiom R$^{++}$:} if $\Sigma_\alpha (\alpha <\omega_1)$ are
stationary subsets of $[X]^{<\omega_1}$ and
$\Gamma \subseteq [X]^{<\omega_2}$
is tight and cofinal, then there is a $Y\in \Gamma$ such
that $\cP (Y)\cap \Sigma_\alpha$ is stationary in
$[Y]^{<\omega_1}$ for each $\alpha <\omega_1$.
\endproclaim
\medskip

Fleissner introduced Axiom R in \cite{Fl} and showed it held
in the usual model for PFA.
\bigskip

\noindent {\bf $\pmb{\Sigma^+}$:} {\it
Suppose $X$ is a countably tight compact space,
$\cL = \{ L_\alpha \}_{\alpha < \omega_1} $ a collection
of disjoint compact sets such that each $L_\alpha$ has a
neighborhood that meets only countably many $L_\beta$'s, and
$\cV$ is a family of $\le \aleph_1$ open subsets of $X$ such that:
\roster
\item "{\it a)}" $\bigcup \cL \subseteq \bigcup \cV $
\item "{\it b)}" For every $V \in \cV$ there is an open $U_V$ such that
$\oV \subseteq U_V $ and $U_V$ meets only countably many members
of $\cL$.
\endroster
\medskip

Then $\cL = \bigcup\limits_{n<\omega} \cL_n$, where each $\cL_n$
is a discrete collection in $\bigcup \cV$. }
\medskip

Balogh \cite{B$_1$} proved that MA$_{\omega_1}$ implies the
restricted version of $\Sigma^+$ in which we take the $L_\alpha$'s
to be points.  We will call that ``$\Sigma'$".
\medskip

\proclaim{Definition}
A space is (strongly) $\kappa$-collectionwise Hausdorff
if for each closed discrete subspace $\{ x_d \}_{d\in D}$,
$|D|\le \kappa$, there is a disjoint (discrete) family of open
sets $\{ U_d \}_{d\in D }$ with $x_d \in U_d$.  A space is
(strongly) collectionwise Hausdorff if it is (strongly)
$\kappa$-collectionwise Hausdorff for all $\kappa$.
\endproclaim
\medskip

It is easy to see that normal $(\kappa-)$~collectionwise
Hausdorff spaces are strongly $(\kappa-)$~collectionwise
Hausdorff.
\medskip

Balogh \cite{B$_2$} proved:
\medskip

\proclaim{Lemma 2}
MA$_{\omega_1}$ + Axiom R implies locally compact hereditarily
strongly \break
$\aleph_1$-collectionwise Hausdorff spaces which do not
include a perfect pre-image of $\omega_1$ are paracompact.
\endproclaim
\medskip

The consequences of MA$_{\omega_1}$ he used are $\Sigma'$ and
Szentmikl\'ossy's result \cite{S} that {\it compact spaces with no
uncountable discrete subspaces are hereditarily Lindel\"of}. Our
plan is to find a model in which these two consequences and Axiom R
hold, as well as normality implying (strongly)
$\aleph_1$-collectionwise Hausdorffness for the spaces under
consideration.  The model we will consider is of the same genre as
those in \cite{LTo}, \cite{L$_2$}, \cite{FTT}, \cite{LT}, and \cite{T$_3$}. One starts
off with a particular kind of Souslin tree $S$, a {\it coherent}
one, which is obtainable from $\diamondsuit$ or by adding a Cohen
real.  One then iterates in standard fashion as in establishing
MA$_{\omega_1}$ or PFA, but omitting partial orders that adjoin
uncountable antichains to $S$. In the PFA case for example, this
will establish {\it PFA(S)}, which is like PFA except restricted to
partial orders that don't kill $S$. In fact it will also establish
PFA$^{++}$(S), which is the corresponding modification of
PFA$^{++}$.  We then force with $S$.
For more information on such
models, see \cite{Mi} and \cite{L$_1$}. We use {\it {\rm PFA}$^{++}(S)[S]$ implies $\varphi$} to
mean that whenever we force over a model of PFA$^{++}(S)$ with $S$,
$\varphi$ holds. Similarly for PFA$(S)[S]$, etc.
\medskip

In \cite{T$_3$} it is established that:
\medskip

\proclaim{Lemma 3} PFA(S)[S] implies that locally compact normal
spaces are $\aleph_1$-collectionwise Hausdorff.
\endproclaim
\medskip

By doing some preliminary forcing (as in \cite{LT}), one can
actually get full collectionwise Hausdorffness, but we won't
need that here.
\medskip

We will assume all spaces are Hausdorff, and use
``$X^*$'' to refer to the one-point compactification of a locally
compact space $X$.
\medskip

There is a bit of a gap in Balogh's proof of Lemma 2.
Balogh asserted that:
\medskip

\proclaim{Lemma 4}
If $X$ is locally compact and does not include a perfect
pre-image of $\omega_1$, then $X^*$ is countably tight.
\endproclaim
\medskip

\noindent and referred to \cite{B$_1$} for the proof.
However in \cite{B$_1$}, he only proved this for the case in
which $X$ is countably tight.  It is not obvious that that
hypothesis can be omitted, but in fact it can.
We need a definition and lemma.
\medskip

\proclaim{Definition}
A space $Y$ is $\pmb\omega${\bf -bounded} if each
separable subspace of $Y$ has compact closure.
\endproclaim
\medskip

\proclaim{Lemma 5}
\cite{G}, \cite{Bu}.
If $Y$ is $\omega$-bounded and does not include a perfect
pre-image of $\omega_1$, then $Y$ is compact.
\endproclaim
\medskip

We then can establish Lemma 4 as follows.
\medskip

\remark{Proof}
By Lemma 5, every $\omega$-bounded subspace of $X$ is compact.
By \cite{B$_1$}, it suffices to show $X$ is countably
tight.  Suppose, on the contrary, that there is a $Y\subseteq X$
which is not closed, but is such that for all countable
$Z\subseteq Y$, $\oZ \subseteq Y$.
Since $X$ is a $k$-space, there is a compact $K$ such
that $K \cap Y$ is not closed.  Then $K\cap Y$ is not
$\omega$-bounded, so there is a countable $Z\subseteq K\cap Y$ such
that $\oZ \cap K \cap Y$ is not compact.
But $\oZ \subseteq Y$, so $\oZ \cap K \cap Y = \oZ \cap K$,
which is compact, contradiction.
\medskip

Lemma 3 takes care of the hereditary strong $\aleph_1$-collectionwise
Hausdorffness we need, since if open subspaces are
$\aleph_1$-collectionwise Hausdorff, all subspaces are, and open
subspaces of locally compact spaces are locally compact.
The proposition that
\medskip

\noindent
$\pmb{\Sigma}$: {\it in a compact countably
tight space, locally countable subspaces of size $\aleph_1$
are $\sigma$-discrete.}
\medskip

\noindent 
is implied by PFA$(S)[S]$ was announced by Todorcevic in the
Toronto Set Theory Seminar in 2002.
\endremark
\medskip

From $\Sigma$ it is standard to get the result of Szentmikl\'ossy
quoted earlier: since the compact space has no uncountable discrete
subspace, it has countable tightness. If it were not hereditarily
Lindel\"{o}f, it would have a right-separated subspace of size
$\aleph_{1}$. But $\Sigma$ implies it has an uncountable discrete
subspace, contradiction.
\medskip

$\Sigma'$ is established by a minor variation of the forcing for
$\Sigma$.  A proof exists in the union of \cite{L$_2$} and \cite{FTT}.
$\Sigma^{+}$, however, is not so clear, and has not yet
been proved from PFA$(S)[S]$. Thus, instead of using it to get
$\aleph_{1}$-collectionwise Hausdorffness in locally compact normal
spaces with no perfect pre-image of $\omega_{1}$, as we did in an
earlier version of this paper, we are instead quoting Lemma 3, which
is a new result of the second author.
\medskip

Thus all we have to do is prove that PFA$^{++}$(S)[S] implies Axiom R.
%
In order to prove that PFA$^{++}$(S)[S] implies Axiom R, we first
note that a straightforward argument using the forcing {\it {Coll}}
$(\omega_1, X)$ (whose conditions are countable partial functions
from $\omega_1$ to $X$, ordered by inclusion) shows that
PFA$^{++}$(S) implies Axiom R$^{++}$.
\medskip

It then suffices to prove:
\medskip

\proclaim{Lemma 6} If Axiom R$^{++}$ holds and $S$ is a Souslin
tree, then Axiom R$^{++}$ still holds after forcing with $S$.
\endproclaim

\remark{Proof}
First note that if $X$ is a set, $P$ is a c.c.c. forcing and
$\tau$ is a $P$-name for a tight cofinal subset of
$[X]^{<\omega_2}$, then the set of $a \in [X]^{<\omega_2}$ such
that every condition in $P$ forces that $a$ is in the realization of
$\tau$ is itself tight and cofinal. The tightness of this set
is immediate. To see that it is cofinal, let $b_0$ be any set
in $[X]^{<\omega_2}$. Define sets $b_\alpha\ (\alpha\le \omega_1)$
and $\sigma_\alpha\ (\alpha < \omega_1)$
recursively by letting $\sigma_\alpha$ be a $P$-name for a
member of the realization of $\tau$ containing $b_\alpha$ and letting
$b_{\alpha+1}$ be the set of members of $X$ which are forced by some
condition in $P$ to be in $\sigma_\alpha$. For limit ordinals
$\alpha\le \omega_1$, let $b_\alpha$ be the union of the
$b_\beta\ (\beta<\alpha)$. Then $b_{\omega_1}$ is forced by every
condition in $P$ to be in $\tau$.

Since we are assuming that the Axiom of Choice holds, Axiom R$^{++}$
does not change if we require $X$ to be an ordinal. Fix an ordinal
$\gamma$ and let $\rho_\alpha (\alpha < \omega_1)$ be $S$-names
for stationary subsets of $[\gamma]^{<\omega_1}$. Let $T$ be a tight
cofinal subset of $[\gamma]^{<\omega_2}$. For each countable ordinal
$\alpha$ and each node $s\in S$, let $\tau_{s,\alpha}$ be the set of
countable subsets $a$ of $\gamma$ such that some condition in $S$
extending $s$ forces that $a$ is
in the realization of $\rho_\alpha$. Applying Axiom R$^{++}$, we have
a set $Y\in [\gamma]^{<\omega_2}$ such that each
$\cP (Y)\cap \tau_{s,\alpha}$ is stationary in $[Y]^{<\omega_1}$.

Since $S$ is c.c.c., every club subset of $[Y]^{<\omega_1}$ that
exists after forcing with $S$ includes a club subset
of $[Y]^{<\omega_1}$ existing in the ground model. Letting
$(\rho_{\alpha})_G$ (for each $\alpha < \omega_1)$
be the realization of $\rho_\alpha$, we have by genericity then
that after forcing with $S$, each $\cP (Y)\cap (\rho_{\alpha})_G$
will be stationary in $[Y]^{<\omega_1}$.
\endremark
\medskip

This completes the proof of Theorem 1.
\medskip

We do not know the answer to the following question; a positive
answer would likely enable us to dispense with Axiom R, and possibly
with the supercompact cardinal.
\medskip

\definition{Problem}
Does MA$_{\omega_1}$ imply every locally compact,
hereditarily strongly collectionwise Hausdorff space which
does not include a perfect pre-image of $\omega_1$ is paracompact?
\enddefinition
\medskip

We also do not know whether in our main result, we can replace
``perfect pre-image of $\omega_1$'' by ``copy of $\omega_1$''.
\medskip

\remark{Remark}
That PFA$(S)[S]$ does not imply Axiom R is proved in \cite{T$_2$}.
\endremark
\medskip

The problem of finding in models of PFA$(S)[S]$ necessary and
sufficient conditions for locally compact normal spaces to be
paracompact is studied in \cite{T$_4$} by extending the methods of
\cite{B$_2$} and this note.

%
%
%
%

\bigskip

{\rm Paul Larson, Department of Mathematics, Miami University,
Oxford, Ohio 45056.}

{\it e-mail address:} {\rm larsonpb\@muohio.edu}
\bigskip

{\rm Franklin D. Tall, Department of Mathematics, University of
Toronto, Toronto, Ontario M5S 2E4, CANADA}

{\it e-mail address:} {\rm f.tall\@utoronto.ca}

\end{document}